\documentclass[11pt,leqno]{article} 
\usepackage{amsmath}
\usepackage{amsfonts}
\begin{document}
\baselineskip=1.1\baselineskip 
\newcommand{\Cal}{\cal}  %For compatibility with Plain.
\newcommand{\goth}{\mathfrak}  %For compatibility with Plain.
\newcommand{\ind}{{\mathbf 1}}    
\newcommand{\ds}{\displaystyle}
\newcommand{\sss}{\scriptstyle}
\newcommand{\rr}{I\!\!R}
\newcommand{\Reals}{\mathbb R} 
\newcommand{\Dom}{{\mbox{Dom}}}
\newtheorem{thm}{Theorem}[section]   %Latex2e does not admit \th
\newtheorem{re}[thm]{Remark}
\newtheorem{co}[thm]{Corollary}
\newtheorem{pr}[thm]{Proposition}
\newtheorem{de}[thm]{Definition}
\newtheorem{lm}[thm]{Lemma}
\newtheorem{example}[thm]{Example}
% Modifying spaces in the description environment:
\renewenvironment{description}
  {\begin{list}{}%
%   {\setlength{\itemsep}{0pt}%
%    \setlength{\itemindent}{0pt}%
%    \setlength{\listparindent}{0pt}%
%    \setlength{\leftmargin}{\labelwidth+\labelsep}}}%
   {\itemsep=0pt%
   \itemindent=0pt%
   \listparindent=0pt%
   \labelwidth=20pt%
   \leftmargin=30pt}}%
  {\end{list}}%
% Modifying the default numbering of equations:
\renewcommand{\theequation}{\thesection.\arabic{equation}}
\def\squarebox#1{\hbox to #1{\hfill\vbox to #1{\vfill}}}
\newcommand{\qed}{\hspace*{\fill}
\vbox{\hrule\hbox{\vrule\squarebox{.45em}\vrule}\hrule}\smallskip}
\def \ci {\mathop{\hbox {\vrule height .4pt depth 0pt width 0.5cm
\hskip -0.3cm \vrule height 10pt
depth 0pt \hskip 0.1cm \vrule height 10pt depth 0pt \hskip 0.2cm}}}
\def\E{\mbox{\rm E}}
\def\supp{\mbox{\rm supp}}
\title{Stochastic differential equations
with boundary conditions
\\
driven by a Poisson noise
\vspace{1cm}}
\author{
Aureli Alabert
\footnote{Supported by grants SGR99-87 of CIRIT and BFM2000-0009 of DGESIC}
\vspace{0.2truecm}
\\
Departament de Matem\`atiques \\
Universitat Aut\`onoma de Barcelona \\
08193 Bellaterra, Catalonia \\
e-mail: alabert@mat.uab.es
\\
\and
Miguel \'{A}ngel Marmolejo
\addtocounter{footnote}{-1}
\footnotemark
\addtocounter{footnote}{-1}\hskip0.15cm
\footnote{$^*$Supported by a grant of Universidad del Valle, Cali-Colombia}
\vspace{0.2truecm}
\\
Departamento de Matem\'aticas\\
Universidad del Valle \\
Cali, Colombia \\
e-mail: mimarmol@univalle.edu.co} 
\maketitle
\begin{abstract}
\noindent
We consider one-dimensional stochastic differential equations with a 
boundary condition, driven by a Poisson process. We study existence 
and uniqueness of solutions and the absolute continuity of the law of the
solution. In the case when the coefficients are linear, we give an 
explicit form of the solution and study the reciprocal process property.
\end{abstract}

\vspace{1truecm}


\noindent
  {\bf Short title:} Poisson SDE with boundary conditions
\thispagestyle{empty}
\vfill\eject
\null
\thispagestyle{empty}
\vfill\eject
\setcounter{page}{1}
\section{Introduction}
\label{intro}
  Stochastic differential equations (s.d.e.) with 
  boundary conditions driven by a Wiener process 
  have been extensively studied in the last fifteen
  years, both in the ordinary and the partial differential
  cases. We highlight the papers of Ocone and Pardoux \cite{Ocone89}, 
  Nualart and Pardoux \cite{Nualart91}, 
  Donati-Martin \cite{Donati-Martin92}, Buckdahn and Nualart \cite{Buckdahn93}, and
  Alabert, Ferrante and Nualart \cite{AFN}. These equations arise from 
  the usual ones when we replace the customary initial condition by a functional relation 
  $h(X_0,X_1)=0$ between two variables of the solution process 
  $X$, which is only considered in the bounded 
  time interval $[0,1]$. Features that have been considered include existence
  and uniqueness, 
  absolute continuity of the laws, numerical approximations,
  and Markovian-type properties.
\par
  Recently, in our work \cite {AM}, we have considered boundary value 
  problems where the stochastic integral with respect to the Wiener process 
  is replaced by an additive Poisson 
  perturbation $N_t$ 
$$
\left\{
\begin{array}{l}
{\ds
X_t=X_{0}+\int_{0}^{t}f(r , X_{r})\, d r+N_{t}
\ ,\quad
t\in[0,1]
\ ,
}
\\ [2mm]
 X_0= \psi (X_1)
 \ ,
\end{array}
\right.
$$
  where the boundary condition is written in a
  more manageable form. 
  We established an existence and uniqueness result, 
  studied the absolutely continuity of the laws, and 
  characterised several classes of 
  coefficients $f$ which lead to 
  the so-called reciprocal property of the solution.
\par
  Let us recall that $X=\{X_t,\ t\in[0,1]\}$ is a 
  {\sl reciprocal process\/} if for all times $0\le s<t \le 1$, the 
  families or random variables $\{X_u,\ u\in [s,t]\}$ and
  $\{X_u,\ u\in [0,1]-[s,t]\}$ are conditionally independent 
  given $X_s$ and $X_t$. This property is weaker than 
  the usual Markov property.
\par
  Interest
  in reciprocal processes dates back to Bernstein \cite{Bernstein32} (they are also
  called Bernstein processes by physicists)
  because of their 
  role in the probabilistic interpretation of quantum mechanics.
  It is by far not true
  that all s.d.e.\ with boundary conditions
  give rise to reciprocal processes and it is also false that 
  a general reciprocal process
  could be represented as the solution of some sort of first order boundary value
  problem, no matter which type of driving process is taken. 
  Nevertheless, it is interesting to try to find 
  in which cases the probabilistic dynamic representation given by a
  first order s.d.e., together with a suitable boundary relation,
  is indeed able to represent a reciprocal process.
\par
  In this paper, we develop the same program as in 
  our previous paper \cite{AM}, but with a multiplicative Poisson
  perturbation. Specifically, we consider the equation
\begin{equation}
\label{B}
\left\{
\begin{array}{l}
{\ds
X_t=X_{0}+\int_{0}^{t}f(r , X_{{r}^{-}})\, d r
+\int_{0}^{t}F(r,X_{{r}^{-}})\,dN_{r} 
\ ,\quad
t \in [0,1]
\ ,
}
\\ [2mm]
 X_0= \psi (X_1)
 \ ,
\end{array}
\right.
\end{equation}
  where $f,F\colon[0,1]\times \Reals \rightarrow \Reals$  and
  $\psi\colon\Reals \rightarrow \Reals$ are measurable 
  functions satisfying certain hypotheses, 
  and $N=\{N_t,\ t \geq 0 \}$ is a Poisson
  process with intensity $1$. 
\par
  Due to the boundary condition, the solution will anticipate
  any filtration to which $N$ is adapted, and therefore 
  the stochastic integral appearing in the equation is,
  strictly speaking, an anticipating integral. However,
  the bounded variation character of the Poisson process
  permits to avoid most of the technical difficulties
  of the anticipating stochastic integrals with respect to
  the Wiener process. 
\par
  Equation (\ref{B}) is a ``forward equation''. One can 
  also consider the ``backward equation'' 
\begin{equation}
\label{D}
\left\{
\begin{array}{l}
{\ds
X_t=X_{0}+\int_{0}^{t}f(r , X_{r})\, d r
+\int_{0}^{t}F(r,X_{r})\,dN_{r} 
\ ,\quad
t \in [0,1]
\ ,
}
\\ [2mm]
 X_0= \psi (X_1)
 \ ,
\end{array}
\right.
\end{equation}
and the Skorohod-type equation
\begin{equation}
\label{C}
\left\{
\begin{array}{l}
{\ds
X_t=X_{0}+\int_{0}^{t}b(r , X_{r})\, d r
+\int_{0}^{t}B(r,X_{r})\,\delta \tilde N_{r} 
\ ,\quad
t \in [0,1]
\ ,
}
\\ [2mm]
 X_0= \psi (X_1)
 \ ,
\end{array}
\right.
\end{equation}
where $\delta \tilde N_{r}$ denotes the Skorohod integral with 
respect to the compensated Poisson process.
  While the stochastic integrals in 
  (\ref{B}) and (\ref{D}) are no more than Stieljes integrals,
  the Skorohod integral operator is defined by means of the 
  chaos decomposition on the canonical Poisson space.
  We refer
  the reader to \cite{LSV}, \cite{NV} or \cite{NV2} for 
  an introduction to the canonical Poisson space, the chaos
  decomposition and the Skorohod integral.
\par
The paper is organised as follows. Section \ref{flows} 
is devoted to the study of the stochastic flow (initial condition problem) 
associated with 
s.d.e.\ (\ref{B}), which will give us the preliminary results 
  needed for the boundary condition case. 
In Section \ref{eqwbc} we study existence, uniqueness, regularity
and absolute continuity of the solution to the problem
(\ref{B}).
In both Sections \ref{flows} and \ref{eqwbc}, the case of linear
equations is studied as a special example.
In Section \ref{reciprocal}, we find some sufficient conditions for the 
solution of the linear equation to enjoy the reciprocal property.
In the final Section \ref{backsko}, the relation of the forward equation
(\ref{B}) with the backward equation
(\ref{D}) and the Skorohod equation (\ref{C}) is established.
The linear equation is again considered with particular attention, and the 
chaos decomposition of the solution is computed in two very 
simple special cases.
\par
\bigskip
  We will use the notation $\partial_i g$ for the derivative
  of a function $g$ with respect to the $i$-th coordinate, $g(s^-)$ 
  and $g(s^+)$ for $\lim_{t\uparrow s}g(t)$ and $\lim_{t\downarrow s}g(t)$
  respectively, and the acronym {\sl c\`adl\`ag\/} for 
  ``right continuous 
  with left limits''.
  Throughout the paper, we employ the usual convention that
  a summation and a product with an empty set of indices 
  are equal to zero and one, respectively. 

\setcounter{equation}{0}
\section{Stochastic flows induced by Poisson equations}
\label{flows}
  Let $N=\{N_t,\ t\ge 0\}$ be a standard Poisson process with intensity 1 defined
  on some probability space $(\Omega,\goth F,P)$; that means, 
  $N$ has independent increments, $N_t-N_s$ has a Poisson law with parameter
  $t-s$, 
  $N_0\equiv 0$, and all its paths are integer-valued, non-decreasing, c\`adl\`ag, 
  with jumps of size 1.
\par
  Throughout the paper, $S_n$ will denote the jump times of $N$:
$$
  S_n(\omega):=\inf\{t\ge 0:\ N_t(\omega)\ge n\}
  \ .
$$
  The sequence $S_n$ is strictly increasing to infinity, and  
  $\{N_t=n\}=\{S_n\le t<S_{n+1}\}$.
\par
\bigskip
  Let us consider the pathwise equation
\begin{equation}
\label{FB}
\varphi_{st}(x)= x+
\int_{s}^{t}f(r,\varphi_{s{r}^{-}}(x))\,dr +
\int_{s}^{t}F(r, \varphi_{s{r}^{-}}(x))\,dN_{r}
\ ,
\quad 0\le s\le t\le 1
\ ,
\end{equation}
  where $x\in\Reals$, 
  and assume that $f,F\colon[0,1]\times\Reals\rightarrow\Reals$
  are measurable functions such that $f$ satisfies
\begin{description}
\item [$(H_1)$]
  $\exists K_1>0:\ \forall t\in[0,1],\ \forall x,y\in\Reals,\
  |f(t,x)-f(t,y)|\le K_1|x-y|$,
\item [$(H_2)$]
  $M_1:=\sup\limits_{t\in[0,1]}|f(t,0)|<\infty$.
\end{description}
\noindent
  For every $x\in\Reals$, denote by $\Phi(s,t;x)$ the solution to the 
  deterministic equation
\begin{equation}
\label{flowdet}
  \Phi(s,t;x)=x+\int_s^t f(r,\Phi(s,r;x))\,dr
  \ ,\quad 0\le s\le t\le 1
  \ .
\end{equation}
All conclusions of the following lemma 
are well known or easy to show:
\begin{lm}
\label{prflowdet}
  Under hypotheses $(H_1)$ and $(H_2)$,
  there exists a unique solution $\Phi(s,t;x)$ of equation
  (\ref{flowdet}). Moreover:
\begin{description}
\item [$1)$]
  For every $0\le s\le t\le 1$, and every $x\in\Reals$,
$
  |\Phi(s,t;x)|\le (|x|+M_1)e^{K_1(t-s)}
  .
$
\item [$2)$]
  For every $0\le s\le r\le t\le 1$, and every $x\in\Reals$,
$
  \Phi(r,t;\Phi(s,r;x))=\Phi(s,t;x)
  .
$
\item [$3)$]
  For every $0\le s\le t\le 1$, and every $x_1,x_2\in\Reals$ with
  $x_1<x_2$,
$$
  (x_2-x_1)e^{-K_1(t-s)}
  \le
  \Phi(s,t;x_2)-
  \Phi(s,t;x_1)
  \le
  (x_2-x_1)e^{K_1(t-s)}
  \ .
$$
  In particular, for every $s,t$, the function $x\mapsto\Phi(s,t;x)$
  is a homeomorphism from $\Reals$ into $\Reals$.
\item [$4)$]
  If $G\colon[0,1]\times\Reals\rightarrow\Reals$ has continuous
  partial derivatives, then for every $0\le s\le t\le 1$,
\begin{align*} 
  G(t,\Phi(s,t;x))=
  G(s,x)
  +\int_s^t
  \Big[
  &
  \partial_1 G(r,\Phi(s,r;x))
  %\cr
  %&
  +\partial_2 G(r,\Phi(s,r;x))
  f(r,\Phi(s,r;x))\Big]
  \,dr
  \ .
  \quad \qed
\end{align*}
\end{description}
\end{lm}
\par
\bigskip
  Using Lemma \ref{prflowdet} one can prove easily the following
  analogous properties for equation (\ref{FB}):
\begin{pr}
\label{generalflujo}
Assume that  
$f$ satisfies hypotheses $(H_1)$ and $(H_2)$
with constants $K_1,M_1$.
Then, for each $x \in \Reals$,
there exists a unique 
%c\`{a}dl\`{a}g  
process
$\varphi(x)=\{\varphi_{st}(x),\ 0\le s\le t\le 1\}$
that solves (\ref{FB}).
Moreover:
\begin{description}
\item [$(1)$]
  If $F$ satisfies hypotheses $(H_1)$ and $(H_2)$ with constants
  $K_2$ and $M_2$, then
  for every $0\leq s \leq t \leq 1$ and every $x \in \Reals$:
$$
|\varphi_{st}(x)|\leq 
\Big[ |x|+(M_1+M_2)(N_t-N_s+1) \Big]
(1+K_2)^{(N_t-N_s)}e^{K_1}
\ .
$$
\item [$(2)$]
For every $0\leq s \leq r \leq t \leq 1$ and every 
$x \in \Reals$,
$
\varphi_{rt}(\varphi_{sr}(x))=
\varphi_{st}(x).
$
\item [$(3)$]
  If there exist constants $-1\le k_2\le K_2$ such that
$$
  k_2(x-y)\le F(t,x)-F(t,y)\le K_2(x-y)
  \ ,\quad
  t\in[0,1]
  \ ,\quad
  x>y
  \ ,
$$
  then for all $0\leq s \leq t \leq1$, 
  and all $x_1, x_2 \in \Reals$  with $x_1<x_2$,
$$
(1+k_2)^{N_t-N_s}e^{-K_1(t-s)} \leq
\frac{\varphi_{st}(x_2)-\varphi_{st}(x_1)}
{x_2-x_1}
\leq
(1+K_2)^{N_t-N_s}e^{K_1(t-s)}
\ ,
$$
  with the convention $0^0=1$.
  In particular, if $k_2>-1$, then for each 
$0\leq s \leq t \leq 1$ the function 
$x \mapsto \varphi_{st}(x)$ is a random homeomorphism from $\Reals$ 
into $\Reals$.
\item [$(4)$]
  Suppose that 
  $G:[0,1]\times \Reals \rightarrow \Reals$ has 
  continuous partial derivatives.
  Then, for all
$0 \leq s \leq t \leq 1 $,
\begin{align*}
G(t,\varphi_{st}(x))=G(s,x)+
&
\int_{s}^{t} \Big[ 
\partial_{1}G(r,\varphi_{s {r}^-}(x))
+\partial_{2}G(r,\varphi_{s {r}^-} (x))
f(r, \varphi_{s {r}^-}(x)) \Big]\,dr
\\
+&\int_{{s}^{+}}^{t} \Big[ G(r,\varphi_{s r^-}(x)+F(r,\varphi_{sr^-}(x)))-
G(r, \varphi_{s{r}^{-}}(x))\Big]\,dN_{r}\ .
\quad \qed
\end{align*}
\end{description}
\end{pr}
  By solving equation (\ref{flowdet}) between jumps, 
  the value $\varphi_{st}(\omega,x)$ can be found recursively in terms
  of $\Phi$: If $s_1=S_1(\omega),\dots, s_n=S_n(\omega)$ are the jump
  times of the path $N(\omega)$ on $(s,1]$, then
\begin{align}
\label{phiPhi}
  \nonumber
  \varphi_{st}(x)
  ={}
  &
  \Phi (s,t;x)\ind_{[s,s_1)}(t)
  +\sum_{i=1}^{n-1}
  \Phi (s_i,t;\varphi_{ss_i^-}(x)+F(s_i,\varphi_{ss_i^-}(x)))
  \ind_{[s_i,s_{i+1})}(t)
  \\
  &
  +\Phi (s_n,t;\varphi_{ss_n^{-}}(x)+F(s_n,\varphi_{ss_n^-}(x)))\ind_{[s_n,1]}(t)
  \ .
\end{align}
  Notice that the paths $t\mapsto\varphi_{st}(x)$ ($t\ge s$) are c\`adl\`ag and 
  $\varphi_{st}(x)-\varphi_{st^{-}}(x)=F(t,\varphi_{st^-})(N_t-N_t^-)$.
\begin{example}
\label{l}
(Linear equation).
{\rm
Let $f_1,f_2,F_1,F_2\colon [0,1]\rightarrow \Reals$ be
continuous functions, and $x \in \Reals$. 
Consider equation (\ref{FB}) with $s=0$ and linear coefficients:
\begin{equation}
\label{FBlineal}
\varphi_{t}(x)=  x+\int_{0}^{t}
[f_1(r)+f_2(r)\varphi_{{r}^{-}}(x) ]\,dr
+\int_{0}^{t}[F_1(r)+F_2(r)\varphi_{{r}^{-}}(x)]
\,dN_{r}\ , \quad 0\le t\le 1
\ .
\end{equation}
We can describe the solution of this equation as
follows: Set $S_0:=0$ and let
$0<S_1<S_2<...$ be the jumps of Poisson process.
For $t \in [S_{i},S_{i+1})$, $i=0,1,2,\dots$,
$$
\varphi_{t}(x)=\varphi_{S_i}(x)+\int_{S_i}^{t}
[f_{1}(r)+f_{2}(r)\varphi_{{r}^{-}}(x) ]\,dr \ .
$$
Applying Proposition \ref{generalflujo}(4)
with $G(t,x)=A(t)^{-1}x$, where
$$
A(t)=\exp \Big\{ \int_{0}^{t}f_2(r)\,dr \Big\}\ ,
$$
we obtain
\begin{equation}
\label{phit}
\frac {\varphi_{t}(x)}{A(t)}=
\frac {\varphi_{S_{i}}(x)}{A(S_i)}+
\int_{S_{i}}^{t}\frac {f_1(r)}{A(r)}\,dr\ .
\end{equation}
On the other hand, for $i=1,2,3,\dots$,
\begin{align}
\nonumber
\varphi_{S_i}(x)
&=\varphi_{S_i^-}(x)+F_1(S_i)+
F_2(S_i)\varphi_{S_i^-}(x) 
\\
\label{phis}
&=F_1(S_i)+[1+F_2(S_i)]\varphi_{S_i^-}(x)\ .
\end{align}
From (\ref{phit}) and (\ref{phis}), it follows that, for $t \in [0,S_1)$,
$$
\frac {\varphi_t (x)}{A(t)}=x+\int_{0}^{t}
\frac {f_1(r)}{A(r)}\,dr\ ,
$$
and that for $t \in [S_i,S_{i+1})$, $i=1,2,\dots$,
\begin{align*}
\frac {\varphi_{t}(x)}{A(t)}=
&\Big[ x+\int_{0}^{S_1}
\frac {f_1(r)}{A(r)}\,dr\Big]\prod_{j=1}^{i}
(1+F_2(S_j)) 
+\Big[\frac {F_1(S_1)}{A(S_1)}+\int_{S_1}^{S_2}
\frac {f_1(r)}{A(r)}\,dr \Big]\prod_{j=2}^{i}
(1+F_2(S_j))+\cdots
\\
&\cdots+\Big[\frac {F_1(S_{i-1})}{A(S_{i-1})}+\int_{S_{i-1}}^{S_i}
\frac {f_1(r)}{A(r)}\,dr \Big]\prod_{j=i}^{i}
(1+F_2(S_j))
+\Big[ \frac {F_1(S_i)}{A(S_i)}+\int_{S_i}^{t}
\frac {f_1(r)}{A(r)}\,dr \Big]\ .
\end{align*}

When $F_2(t)\neq -1$ for almost all $t \in [0,1]$ with respect to 
Lebesgue measure, we can also write 
the solution as follows:
$$
\varphi_{t}(x)=\eta_{t} \Big[x+\int_{0}^{t}
\frac {f_1(r)}{\eta_{r}}\,dr+\int_{0}^{t}
\frac {F_1(r)}{\eta_{r}}\,dN_{r} \Big]
\ ,
\quad\mbox{a.s.,}
$$                                
where
$$
\eta_{t}=A(t)\prod_{0<S_i\leq t}[1+F_2(S_i)]
\ .
\quad\qed
$$
}
\end{example}
\par
Under differentiability assumptions on 
$f$ and $F$ we obtain
differentiability properties of the solution
to (\ref{FB}): 
\begin{pr}
\label{derflujogeneral}
Assume that $f$ satisfies the stronger hypotheses:
\begin{description}
\item [$(H_1')$]
$f$, $\partial_{2}f$ are continuous functions.
\item [$(H_2')$]
$\exists K>0:\ 
|\partial_{2}f|\leq K$.
\end{description}
Then 
\begin{description}
\item [$(1)$]
  For every $\omega\in\Omega$ and every $x\in\Reals$, the function
  $t\mapsto \varphi_{st}(\omega,x)$ is differentiable on
  $[s,1]-\{s_1,s_2,\dots\}$,
  where $s_1,s_2,\dots$ are the jump times of $N(\omega)$ on $(s,1]$,
  and
$$
  \frac{d\varphi_{st}(x)}{dt}=f(t,\varphi_{st}(x))
  \ .
$$
\item [$(2)$]
  If moreover $F$ and $\partial_2 F$ are continuous functions, 
  then for every $\omega\in\Omega$ and every $0\le s\le t\le 1$,
  the function $x\mapsto \varphi_{st}(\omega,x)$ is 
  continuously differentiable and
$$
\frac {d\varphi_{st}(\omega,x)}{dx}=
\exp \Big\{ \int_{s}^{t}
\partial_{2}f(r,\varphi_{sr}(\omega,x))\,dr \Big\}
\prod_{s<s_i\leq t}
[ 1+\partial_{2}F(s_i,\varphi_{ss_i^-}(\omega,x))]
\ .
$$
  In particular, when $\partial_2 F>-1$,
  $x \mapsto \varphi_{st}(x)$ is a
  random diffeomorphism from $\Reals$ into $ \Reals$.

\item [$(3)$]
  Assume moreover that $F$, $\partial_1 F$ and $\partial_2 F$
  are continuous functions.
Fix $0\leq s<t \leq 1$ and $n \in \{1,2,\dots\}$. On the 
set $\{N_t-N_s=n\}$, the mapping $\omega\mapsto \varphi_{st}(\omega,x)$ 
regarded as a
function $\varphi_{st}(s_1,s_2,\dots,s_n;x)$ defined
on 
$\{s<s_1<s_2<\cdots<s_n\le t\}$ (where $s_j=S_j(\omega)$ are the jump 
times of
$N(\omega)$ on $(s,t]$), 
is continuously differentiable and, for every 
$j\in\{1,\dots,n\}$,
\begin{align*}
\frac {\partial \varphi_{st}(x)}{\partial s_j}=
&
\exp \Big\{ \int_{s_j}^{t}
\partial_{2}f(r,\varphi_{sr}(x))\,dr \Big\}
\prod_{i=j+1}^{n}
[1+\partial_{2}F(s_i,\varphi_{ss_i^-}(x))]
\\
&
\times\Big[ -f(s_j,\varphi_{ss_j}(x))+
\partial_{1}F(s_j,\varphi_{ss_j^-}(x))
+
f(s_j,\varphi_{ss_j^-}(x))
[1+\partial_{2}F(s_j,\varphi_{ss_j^-}(x))]\Big]
\ .
\end{align*}
\par
\end{description}
\end{pr}
{\em Proof}:
  It is easy to see for the solution $\Phi$ of (\ref{flowdet}) that 
\begin{align*}
  \partial_{1} \Phi (s,t;x)
  &=
  -f(s,x)\exp \Big\{
  \int_{s}^{t}\partial_{2}f(r,\Phi (s,r;x))\,dr\Big\}
  \ ,
  \\
  \partial_{2}\Phi (s,t;x)
  &=
  f(t,\Phi (s,t;x))
  \ ,
  \\
  \partial_{3}\Phi (s,t;x)
  &=
  \exp \Big\{\int_{s}^{t}
  \partial_{2}f(r,\Phi (s,r;x))\,dr\}
  \ ,  
\end{align*}
  and that these derivatives are continuous on $\{0\le s\le t\le 1\}\times\Reals$.
  Claims (1) and (2) follow from here and representation (\ref{phiPhi}).
\par
  The existence and regularity of the function $\varphi_{st}(s_1,\dots,s_n;x)$ of 
  (3) are also clear from (\ref{phiPhi}). We compute now its derivative 
  with respect to $s_j$. For $n=1$, we get
\begin{align*}
  \frac {d\varphi_{st}(x)}{ds_1}
  =
  {}&
  \partial_{1}\Phi (s_1,t;\varphi_{ss_1}(x))+
  \partial_{3}\Phi (s_1,t;\varphi_{ss_1}(x))
  \frac {d\varphi_{ss_1}(x)}{ds_1}
  \\
  =
  {}&
  \exp \Big\{\int_{s_1}^{t}
  \partial_{2}f(r,\varphi_{sr}(x))\,dr\Big\}
  \\
  {}&
  \times\Big[- f(s_1,\varphi_{ss_1}(x)) +
  \partial_{1}F(s_1,\varphi_{ss_1^-}(x)) 
  +f(s_1,\varphi_{ss_1^{-}}(x))
  [1+\partial_{2}F(s_1,\varphi_{ss_1^{-}}(x))] \Big] 
  \ .
\end{align*}
Suppose that $(3)$ holds for $n=k$. Then, for
$n=k+1$ and $j=1,\dots,k$,
\begin{align*}
\frac {\partial \varphi_{st}(x)}{\partial s_j}
={}&
\partial_{3} \Phi (s_{k+1},t;\varphi_{ss_{k+1}}(x))
\frac {\partial \varphi_{ss_{k+1}}(x)}
{\partial s_j}
\\
={}&\exp \Big\{\int_{s_j}^{t}
\partial_{2}f(r,\varphi_{sr}(x))\,dr \Big\}
\prod_{i=j+1}^{k+1}
[1+\partial_{2}F(s_i,\varphi_{ss_i^-}(x)) ]
\\
{}&\times\Big[ -f(s_j,\varphi_{ss_j}(x))+
\partial_{1}F(s_j,\varphi_{ss_j^-}(x)) 
+f(s_j,\varphi_{ss_j^-}(x))
[1+\partial_{2}F(s_j,\varphi_{ss_j^-}(x))] \Big]
  \ . 
\end{align*}
Taking into account that 
\begin{align*}
\varphi_{ss_{k+1}}(x)
&= \varphi_{ss_{k+1}^-}(x)+
F(s_{k+1},\varphi_{ss_{k+1}^-}(x))
\ ,
\\
\frac {\partial \varphi_{s{s_{k+1}}}(x)}
{\partial s_{k+1}}
&= 
[1+\partial_{2}F(s_{k+1},\varphi_{ss_{k+1}^-}(x)) ]
f(s_{k+1},\varphi_{ss_{k+1}^-}(x))
+\partial_{1}F(s_{k+1},\varphi_{ss_{k+1}^-}(x))
\ ,
\end{align*}
we obtain, for $j=k+1$,
\begin{align*}
\frac {\partial \varphi_{st}(x)}{\partial s_{k+1}}
={}&\partial_{1}\Phi (s_{k+1},t;\varphi_{ss_{k+1}}(x))+
\partial_{3}\Phi (s_{k+1},t;\varphi_{ss_{k+1}}(x))
\frac {\partial \varphi_{ss_{k+1}}(x)}{\partial s_{k+1}} 
\\
={}&\exp \Big \{\int_{s_{k+1}}^{t}
\partial_{2}f(r,\varphi_{sr}(x))\,dr \Big\}
\\
&\times\Big[-f(s_{k+1},\varphi_{ss_{k+1}}(x))+
\partial_{1}F(s_{k+1},\varphi_{ss_{k+1}^-}(x))
\\
&
\phantom{{}\times\Big[}
+f(s_{k+1},\varphi_{ss_{k+1}^-}(x))
[1+\partial_{2}F(s_{k+1},\varphi_{ss_{k+1}^-}(x))]
\Big]
\ . 
\quad\qed   
\end{align*}
\par
In the next proposition we find that under the regularity 
hypotheses of Proposition \ref{derflujogeneral} and 
an additional condition relating $f$ and $F$, the law 
of $\varphi_t(x)$ is a weighted sum of a Dirac-$\delta$
and an absolutely continuous probability.
\begin{pr}
\label{leyflujogeneral}
Let $f$ satisfy hypotheses $(H_1')$ and $(H_2')$ of 
Proposition \ref{derflujogeneral}, and 
assume that $F$, $\partial_1 F$ and $\partial_2 F$ are 
continuous functions. 
Assume moreover that 
\begin{equation}
\label{conditionP}
|f(t,x+F(t,x))-
 f(t,x)[1+\partial_{2}F(t,x)]-
\partial_{1}F(t,x)|> 0
\ ,\quad
\forall t\in[0,1],\ \forall x\in\Reals
\ .
\end{equation}
Let $\varphi(x)=\{\varphi_t(x),\ t\in[0,1]\}$ be the solution
to (\ref{FB}) for $s=0$. Then,
for all $t>0$, the distribution 
function $L$ of $\varphi_{t}(x)$ can be 
written as
$$
L(y)=
e^{-t}L^{D}(y) +
(1-e^{-t})L^{C}(y)
\ ,
$$
with
$$
L^{D}(y)=\ind_{[\Phi (0,t;x),\infty)}(y)\ ,
$$
and
$$
L^{C}(y)=(e^{t}-1)^{-1}
\int_{-\infty}^{y} \sum_{n=1}^{\infty}
\frac {t^{n}}{n!}
h_n(r)\,dr
\ ,
$$
where $h_n$ is the density function
of the law of $\varphi_{t}(x)$ conditioned to $N_t=n$.
\end{pr}
{\em Proof}:
  Let $S_1,S_2,\dots$ be the jump times
  of $\{N_t,\ t\in[0,1]\}$.
  From Proposition \ref{derflujogeneral}(3), on the set
  $\{N_t=n\}$ ($n=1,2,\dots$) we have 
  $
  \varphi_t(x)=G(S_1,\dots,S_n)
  $
  for some continuously differentiable function $G$, and that
\begin{align*}
\partial_{n}G(s_1,...,s_n)
=&\exp \Big\{ \int_{s_n}^{t}
\partial_{2}f(r,\varphi_{r}(x))\,dr \Big\}
\\
&\times\Big[ -f(s_n,\varphi_{s_n}(x))+\partial_{1}
F(s_n,\varphi_{s_n^-}(x))+
f(s_n,\varphi_{s_n^-}(x))[1+\partial_{2}
F(s_n,\varphi_{s_n^-}(x))] \Big]
\ .
\end{align*}
Using
$
\varphi_{s_n}(x)=\varphi_{{s_n^-}}(x)+
F(s_n,\varphi_{{s_n^-}}(x))
$
and condition (\ref{conditionP}), we obtain
$
|\partial_{n}G|>0.
$
\par
  It is known that, conditionally to $\{N_t=n\}$,
  $(S_1,\dots,S_n)$ follows the uniform distribution
  on $D_n=\{0<s_1<\cdots<s_n< t\}$.
  If we define $T(s_1,\dots,s_n)=(z_1,\dots,z_n)$, with
  $z_i=s_i$, $1\le i\le n-1$, and $z_n=G(s_1,\dots,s_n)$, then
  $(Z_1,\dots,Z_n)=T(S_1,\dots,S_n)$ is a random vector with density
$$
  h(z_1,\dots,z_n)=
  n!\,t^{-n}\big|\partial_n s_n
  (z_1,\dots,z_n)\big|
  \ind_{T(D_n)}(z_1,\dots,z_n)
  \ ,
$$
  and therefore $\varphi_t(x)$ is absolutely continuous on $\{N_t=n\}$,
  for every $n\ge 1$,
  with conditional density
$$
  h_n(y)=\ind_{G(D_n)}(y)
  \int\int\cdots\int
  h(z_1,\dots,z_{n-1},y)  \,dz_1\,\dots\,dz_{n-1}
  \ .
$$
  Now,
\begin{align*}
  L(y)
  &
  =
  \sum_{n=0}^{\infty}
  P\big\{\raise2pt\hbox{$\varphi_t(x)\le y$}/
  \raise-2pt\hbox{$N_t=n$}\big\}
  P\{N_t=n\}
  \\
  &
  =
  e^{-t}
  P\big\{\raise2pt\hbox{$\varphi_t(x)\le y$}/
  \raise-2pt\hbox{$N_t=0$}\big\}
  +
  e^{-t}
  \sum_{n=1}^{\infty} \frac{t^n}{n!}
  \int_{-\infty}^y
  h_n(r)\,dr
  \\
  &
  =
  e^{-t}
  \ind_{[\Phi(0,t;x),\infty)}(y)
  +
  e^{-t}
  \int_{-\infty}^y
  \sum_{n=1}^{\infty} \frac{t^n}{n!}
  h_n(r)
  \,dr
  \ ,
\end{align*}
  and the result follows.
\begin{re}
\label{o8}
{\rm
When $f(t,x)\equiv f(x)$, $F(t,x)\equiv F(x)$ and
$f''$ is continuous, condition (\ref{conditionP})
is satisfied if
$$
|f'F-fF'|>
\frac {1}{2} \|f''\|_{\infty}
             \|F\|_{\infty}^{2}\ ,
$$ 
which is the hypothesis used by Carlen and Pardoux in
\cite {Carlen90} (Theorem 4.3) to prove that, in the autonomous 
case, the law of $\varphi_{1}(x)$ is 
absolutely continuous 
on the set
$\{N_1 \geq 1 \}$. 
}
\end{re}
\setcounter{equation}{0}
\section{Equations with boundary conditions}
\label{eqwbc}
  In this section 
  we establish first an easy existence and uniqueness theorem, based
  on Proposition \ref{generalflujo} above, 
  when the initial condition is replaced by a boundary 
  condition.
  Then we prove in this situation the analogue 
  of Propositions \ref{derflujogeneral}(3) and \ref{leyflujogeneral} on the differentiability 
  with respect to the jump times and the absolute continuity of the laws
  (Proposition \ref{derspgeneral} and Theorem \ref{laws2} below, respectively).
\begin{thm}
\label{spB}
Let $f,F\colon[0,1]\times \Reals \rightarrow \Reals$ be
measurable functions such that $f$ 
satisfies hypotheses $(H_1)$ and $(H_2)$ of Section \ref{flows}, with constants $K_1$ and
$M_1$ respectively, and there exists a constant  
$k_2\ge -1$ such that $F(t,x)-F(t,y)\ge k_2(x-y)$, $\forall t\in[0,1]$, $x>y$.
Assume that $\psi\colon \Reals\to\Reals$
satisfies
\begin{description}
\item [$(H_3)$]
$\psi$ is a continuous and non-increasing function.
\end{description}
Then 
\begin{equation}
\label{B2}
\left\{
\begin{array}{l}
{\ds
X_t=X_{0}+\int_{0}^{t}f(r , X_{{r}^{-}})\, d r
+\int_{0}^{t}F(r,X_{{r}^{-}})\,dN_{r} 
\ ,\quad
t \in [0,1]
\ ,
}
\\ [2mm]
 X_0= \psi (X_1)
 \ ,
\end{array}
\right.
\end{equation}
admits a unique solution
$X$, which is a c\`{a}dl\`{a}g process.
\end{thm}
{\em Proof}:
By Proposition \ref{generalflujo}, for
each $x \in \Reals$ there exists a unique c\`{a}dl\`{a}g process
$\varphi (x)=\{\varphi_{t}(x),\ t \in [0,1]\}$
that satisfies the equation
$$
\varphi_{t}(x)=x+\int_{0}^{t}
f(r, \varphi_{{r}^{-}}(x))\,dr+\int_{0}^{t}
F(r,\varphi_{{r}^{-}}(x))\,dN_{r}\ ,
\quad t \in [0,1] 
\ .
$$
From part (3) of the same proposition, for each 
$\omega \in \Omega$, the function
$x \mapsto \varphi_{1}(\omega,x)$ is non-decreasing.
Thus, by hypothesis $(H_3)$,
the function $x \mapsto x-\psi (\varphi_{1}(\omega,x))$
has a unique fixed point, that we define as
$X_0(\omega)$. It follows that $(\ref{B2})$
has the unique solution 
$X_t(\omega)=\varphi_{t}(\omega,X_0(\omega))$.\qed
\begin{re}
%\label{}
{\rm
In general, the condition $k_2\ge -1$ cannot be 
relaxed. For instance, the  problem
$$
\left \{
\begin{array}{l}
\ds
X_t=X_0+\int_{0}^{t} X_{{r}^{-}}\,dr +\int_{0}^{t}
-2X_{{r}^{-}}\,dN_{r} \ ,
\vspace{2truemm}
\\
X_0=1-\frac{1}{e}{X_1}\ ,\quad t \in [0,1]\ ,
\end{array}
\right.
$$
has no solutions. Indeed, the first equality implies
$
X_t=X_0e^{t}(-1)^{N_t}
$ (see Example \ref{l}),
which gives $X_1=-X_0 e$ for 
$N_{1}\in \{1,3,5,\dots\}$, and this is incompatible 
with the boundary condition.
\par
On the other hand, if we change the boundary condition 
to $X_0=\frac {-1}{e}X_1$
the new problem has an infinite number of solutions: 
$$
X_t(\omega)=
\left \{
\begin{array}{ll}
\ds
0\ , &  \mbox{if $N_1(\omega)=0,2,4,\dots$} \cr
x(\omega)e^{t}(-1)^{N_{t}(\omega)}\ , 
& \mbox{if $N_1(\omega)=1,3,5,\dots$}
\end{array}
\right.
$$
where $x (\omega)$ is an arbitrary real number.
\par
Notice that the purpose of $(H_3)$ is to ensure that
$x\mapsto \psi(\varphi_1(x))$ has a unique fix point.
Alternative hypotheses that lead to the same consequence
may be used instead for particular cases. See for instance
the comments at the end of Example \ref{OE}.
\qed
}
\end{re}
\begin{example}
\label{OE}
(Linear equation).
{\rm 
Consider the linear equation
\begin{equation}
\label{LB}
\left \{
\begin{array}{l}
\ds
X_t=X_0+\int_{0}^{t} 
[f_1(r)+f_2(r)X_{{r}^{-}}]\,dr +\int_{0}^{t}
[F_1(r)+F_2(r)X_{{r}^{-}}]\,dN_{r} \ ,\\
X_0=\psi(X_1) \ , \quad t \in [0,1]\ ,
\end{array}
\right.
\end{equation}
where $f_1,f_2,F_1,F_2\colon[0,1] \rightarrow \Reals$ are
continuous functions, and
$\psi \colon\Reals \rightarrow \Reals$ is a continuous and
non-increasing function. Assume
$F_2(t)\geq -1$ for all $t \in [0,1]$. By 
Theorem \ref{spB} there is a unique solution, and using Example \ref{l}
we can describe it as follows:
% which can be
%given explicitly, with the help of Example \ref{l}, as follows:
%Example \ref{l}, this problem has unique solution,
%which can be described as follows:
\par
For $\omega \in \{S_1>1\}$,
$$
X_t=A(t)\Big[ x^{\ast}+\int_{0}^{t}
\frac {f_1(r)}{A(r)}\,dr \Big]
\ ,
$$
where $A(t)=\exp \{\int_{0}^{t}f_2(r)\,dr \}$,
and $x^{\ast}$ solve
$$
x=\psi \Big( A(1)\Big[x+\int_{0}^{1}\frac {f_1(r)}{A(r)}
\,dr \Big]\Big)
\ .
$$
\par
For $\omega \in \{S_n<1<S_{n+1}\}$ $(n\ge1)$ 
and $t \in [S_i,S_{i+1})$
we have
\begin{align*}
\frac {X_{t}(x)}{A(t)}=
&\Big[ X_0+\int_{0}^{S_1}
\frac {f_1(r)}{A(r)}\,dr\Big]\prod_{j=1}^{i}
(1+F_2(S_j)) 
+\Big[\frac {F_1(S_1)}{A(S_1)}+\int_{S_1}^{S_2}
\frac {f_1(r)}{A(r)}\,dr \Big]\prod_{j=2}^{i}
(1+F_2(S_j))+\cdots
\\
&\cdots+\Big[\frac {F_1(S_{i-1})}{A(S_{i-1})}+\int_{S_{i-1}}^{S_i}
\frac {f_1(r)}{A(r)}\,dr \Big]\prod_{j=i}^{i}
(1+F_2(S_j))
+\Big[ \frac {F_1(S_i)}{A(S_i)}+\int_{S_i}^{t}
\frac {f_1(r)}{A(r)}\,dr \Big]\ .
\end{align*}
where $X_0$ solve
$x=\psi (\varphi_1(x))$, with
\begin{align*}
\frac {\varphi_{1}(x)}{A(1)}=
&\Big[ x+\int_{0}^{S_1}
\frac {f_1(r)}{A(r)}\,dr\Big]\prod_{j=1}^{i}
(1+F_2(S_j)) 
+\Big[\frac {F_1(S_1)}{A(S_1)}+\int_{S_1}^{S_2}
\frac {f_1(r)}{A(r)}\,dr \Big]\prod_{j=2}^{i}
(1+F_2(S_j))+\cdots
\\
&\cdots+\Big[\frac {F_1(S_{i-1})}{A(S_{i-1})}+\int_{S_{i-1}}^{S_i}
\frac {f_1(r)}{A(r)}\,dr \Big]\prod_{j=i}^{i}
(1+F_2(S_j))
+\Big[ \frac {F_1(S_i)}{A(S_i)}+\int_{S_i}^{1}
\frac {f_1(r)}{A(r)}\,dr \Big]\ .
\end{align*}

When $F_2(t)> -1$ for almost all $t \in [0,1]$ with respect to 
Lebesgue measure, we can also write 
the solution as follows:
\begin{equation}
\label{explicitLin}
X_{t}(x)=\eta_{t} \Big[X_0+\int_{0}^{t}
\frac {f_1(r)}{\eta_{r}}\,dr+\int_{0}^{t}
\frac {F_1(r)}{\eta_{r}}\,dN_{r} \Big]
\ ,
\quad\mbox{a.s.,}
\end{equation}                                
where
\begin{equation*}
\eta_{t}=A(t)\prod_{0<S_i\leq t}[1+F_2(S_i)]
=
\exp\Big\{\int_0^t f_2(r)\,dr + \int_0^t \log(1+F_2(r))\,dN_r\Big\}
\ .
\quad
\end{equation*}
Finally, we remark that if $-1\le F_2\le 0$, 
the monotonicity condition on $\psi$ can be relaxed to 
$$
x>y\Rightarrow \psi(x)-\psi(y)\le \alpha(x-y)
\ ,
$$ 
with
$\alpha A(1)<1$, because in this case the mapping 
$x\mapsto x-\psi(\varphi_1(\omega,x))$ has still a unique fix point.
\qed
}
\end{example}
Under differentiability assumptions on
$f$, $F$ and $\psi$,
we will obtain differentiability properties of the solution to
$(\ref{B2})$. 
Denote
\begin{align*}
A(s_j,t,X):=&\exp \Big\{ \int_{s_j}^{t}
\partial_{2}f(r,X_{r})\,dr \Big\} 
\prod_{s_j< s_i \leq t }[1+\partial_{2}F(s_i,X_{s_i^-})]
\\
&
\times
\Big[
-f(s_j,X_{s_j})+
f(s_j,X_{s_j^-})[1+\partial_{2}F(s_j,X_{s_j^-})]
\partial_{1}F(s_j,X_{s_j^-})
\Big]\ ,
\end{align*}
and
$$
B(t,X):=\exp \Big\{\int_{0}^{t}
\partial_{2}f(r,X_{r})\,dr \Big\}
\prod_{0< s_i \leq t }[1+\partial_{2}F(s_i,X_{s_i^-})]
\ .
$$
\begin{pr}
\label{derspgeneral} 
Let $f,F\colon[0,1]\times \Reals \rightarrow \Reals$ and
 $\psi\colon\Reals \rightarrow \Reals$ be measurable functions such that
$f$ satisfies hypotheses $(H_1')$, $(H_2')$ of Section \ref{flows}; $F$,
$\partial_1 F$ and $\partial_2 F$ are continuous functions
with $\partial_2 F\geq -1$, and 
\begin{description}
\item [$(H_3')$]
$\psi$ is a continuously differentiable function with $\psi ' \leq 0$.
\end{description}
Let $X=\{X_t,\ t \in [0,1]\}$ be the solution to  
(\ref{B2}). Then:                  
\begin{description}
\item [$(1)$]
Fix $n \in \{1,2,\dots\}$. On the set $\{N_1=n\}$,
$X_0$ can be regarded as a function $X_0(s_1,s_2,...,s_n$), 
defined on $\{0<s_1<\cdots<s_n<1\}$, where
$s_j=S_j(\omega)$ are the jumps of  
$N(\omega)$ in $[0,1]$. 
This function is continuously differentiable, and for 
any $j=1,2,\dots,n$,
$$
\frac {\partial X_0}{\partial s_j}=
\frac {\psi ' (X_1)A(s_j,1,X)}
{1-\psi '(X_1)B(1,X)}
\ .
$$
\item [$(2)$]
Take $t \in (0,1]$ and $n,k \in \{0,1,\dots\}$  such that
$n+k \geq 1$. On the set
$\{N_t=n\}\cap \{N_1-N_t=k\}$, $X_t$ can be regarded as a function 
$X_t(s_1,\dots,s_{n+k})$ 
defined on $\{0<s_1<\cdots<s_{n+k}<1\}$, where
$s_j=S_j(\omega)$ are the jumps of  
$N(\omega)$ in $[0,1]$. 
This function is continuously differentiable, and for 
any $j=1,2,\dots,n+k$,
$$
\frac {\partial X_t}{\partial s_j}=
B(t,X)\frac {\partial X_0}{\partial s_j}
\ind_{\{1,\dots,n+k\}}(j)+
A(s_j,t,X) 
\ind_{\{1,\dots,n\}}(j)
\ .
$$
\end{description}
\end{pr}
{\em Proof}:
Since $X_0=\psi (\varphi_{1}(X_0))$, we have
\begin{align*}
\frac {\partial X_0}{\partial s_j}
&=
\frac {\psi '(\varphi_{1}(X_0))
\frac {\partial \varphi_{1}(x)}{\partial s_j}\Big|_{x=X_0}}
{1-\psi '(\varphi_{1}(X_0))
\frac {d\varphi_{1}(x)}{dx}\Big|_{x=X_0}}
\\
&=
\frac {\psi '(X_1) A(s_j,1,X)}
{1-\psi '(X_1)B(1,X)}
\ .
\end{align*}
On the other hand,
for $X_t=\varphi_{t}(X_0)$,
\begin{align*}
\frac {\partial X_t}{\partial s_j}
=&
\frac {d\varphi_{t}(x)}{dx}\Big|_{x=X_0} 
\frac {\partial X_0}{\partial s_j}  +
\frac {\partial \varphi_{t}(x)}{\partial s_j}
\Big|_{x=X_0}
\\
=&
B(t,X)\frac {\partial X_0}{\partial s_j}
\ind_{\{1,\dots,n+k\}}(j)
+
 A(s_j,t,X)
\ind_{\{1,\dots,n\}}(j)
\ .\quad\qed 
\end{align*}
The following theorem is the counterpart of Proposition \ref{leyflujogeneral}
for the case of boundary conditions. The proof follows the same
lines but using at the end the decomposition
\begin{align*}
  L_{X_{t}}(x)
  ={}&
  P\{X_{t}\leq x,\, N_{1}=0\}
  \\
  &
  +
  \sum_{n=1}^{\infty}
  P\big\{\raise2pt\hbox{$X_{t}\leq x$}/
  \raise-2pt\hbox{$N_{t}=0,\ N_{1}-N_{t}=n$}\big\}e^{-1} 
  \frac{(1-t)^{n} }{n!}
  \\
  &   
  +
  \sum_{n=1}^{\infty}
  P\big\{\raise2pt\hbox{$X_{t}\leq x$}/\raise-2pt\hbox{$N_{t}=n$}\big\}
  e^{-t}
  \frac {t^{n}}{n!}
  \ .
  \quad\qed
\end{align*}
\begin{thm}
\label{laws2} 
Let $f,F\colon[0,1]\times \Reals \rightarrow \Reals$
and $\psi\colon \Reals \rightarrow \Reals  $ satisfy the hypotheses
of Proposition \ref{derspgeneral}. 
Assume in addition that $\psi '<0$ and that condition 
(\ref{conditionP}) holds. 
  Let $x^*$ be the unique solution to
$x=\psi(\Phi(0,1;x))$,
and $X$ the solution to (\ref{B2}).
Then, 
  the distribution function of $X_t$, $t \in (0,1]$, is
$$
  L_{X_t}(x)=e^{-1}L_{X_t}^{D}(x)+(1-e^{-1})L_{X_t}^{C}(x)
  \ ,
$$
  with
$$
L_{X_t}^{D}(x)=\ind_{[\Phi (0,t;x^{\ast}),\infty)}(x)
$$
  and
$$
  L_{X_t}^{C}(x)=
  \frac {e^{-t}}{1-e^{-1}} \Big[
  e^{-(1-t)}\int_{-\infty}^{x} \sum_{n=1}^{\infty}
  \frac {(1-t)^{n}}{n!} h_{0n}(r )\,d r 
  +\int_{-\infty}^{x} \sum_{n=1}^{\infty}
  \frac {h_n(r)}{n!}\, d r \Big]
  \ ,
$$
  where $h_{0n}$ is the density of
  $X_t$ conditioned to $N_{t}=0,N_{1}=n$, and 
  $h_n$ is the density of $X_t$ conditioned
  to $N_t=n$.
  For $t=0$, the formula is also valid taking $h_n\equiv 0$.
\end{thm}
\section {The reciprocal property}
\label{reciprocal}
Let $(\Omega, \goth F,P)$ be a probability space and
let $\goth A_1$, $\goth A_2$ and $\goth B$ 
be sub-$\sigma$-fields of $\goth F$ such that
$
P(A_1\cap A_2|\goth B)=
P(A_1|\goth B)P(A_2|\goth B)
$
for any $A_1 \in \goth A_1$, $A_2 \in \goth A_2$.
Then the $\sigma$-fields $\goth A_1$
and $\goth A_2$ are said to be {\sl conditionally independent with
respect to $\goth B$}.
\begin{de}
{\rm
We say that $X=\{X_t,\ t\in[0,1]\}$ is a {\sl reciprocal process} if for every
$0\leq a<b \leq 1$, the $\sigma$-fields generated by
$\{X_t ,\ t \in [a,b] \}$ and $\{X_t,\ t \in (a,b)^{c} \}$
are conditionally independent with respect to the 
$\sigma$-field generated by $\{X_a,X_b\}$.
}
\end{de}
One can show that if $X$ is a 
Markov process then $X$ is reciprocal, and that the converse
is not true. 
For a proof of this fact, we refer the reader to Alabert and Marmolejo \cite {AM}
(Proposition 4.2),
where we also established the next lemma (Lemma 4.6 of \cite{AM}). 
\par
\begin{lm}
\label{JM}
If $\xi=\{ \xi_t,\ t \in [0,1]\}$ has independent increments
and $g$ is a Borel function, then 
$X:=\{g(\xi_1)+\xi_t,\ t \in [0,1]\}$ is a reciprocal process.
\end{lm}

In our previous work \cite{AM}, we obtained several sufficient
conditions on $f$ for the solution to enjoy the reciprocal 
property when the Poisson noise 
appears additively, namely in 
$$
\left \{
\begin{array}{l}
\ds
X_t=X_0+\int_{0}^{t}f(r,X_{r})\,dr+
\int_{0}^{t}\,dN_{r}\ ,  \\
X_0=\psi (X_1)\ , \quad t \in [0,1]\ .
\end{array}
\right.
$$  
  The main classes of functions $f$ leading to this property
  are those of the form $f(t,x)=f_1(t)+f_2(t)x$ 
  and those which are 1-periodic in the second variable, 
  $f(t,x)=f(t,x+1)$. But we showed with examples that there are
  many more; we also
  obtained conditions on $f$ ensuring that 
  the solution will not be reciprocal.
  In contrast, for equations driven by the Wiener process, 
  conditions which are at the same time
  necessary
  and sufficient 
  have been obtained in a wide variety
  of settings, even with multiplicative noise. 
\par
  With a multiplicative Poisson noise, the techniques currently
  known do not
  seem to allow a general analysis. 
  We will restrict ourselves to linear equations. 
  The main result contained in the next
  theorem is that if both coefficients are truly linear in the
  second variable (i.e. $f(t,x)=f_2(t)x$, and $F(t,x)=F_2(t)x$)
  then the solution is reciprocal. We have not been able to obtain
  necessary conditions, even when considering only the class
  of linear equations. Thus, we have to 
  leave open the study of the general linear case, which for white noise
  driven equations (with boundary conditions also linear) was 
  studied thoroughly in the seminal paper of Ocone and Pardoux \cite{Ocone89}.
\begin{thm}
Let $f_1,f_2,F_1,F_2\colon[0,1]\rightarrow \Reals$ be
continuous functions with $ F_2(t) \geq -1$ for all 
$t \in [0,1]$, and $\psi \colon\Reals \rightarrow \Reals$ 
a continuous and non-increasing function. Let
$X=\{X_t,\ t \in [0,1]\}$ be the solution of  
$$
\left \{
\begin{array}{l}
\ds
X_t=X_0+\int_{0}^{t}
[f_1(r)+f_2(r)X_{{r}^{-}}]\,dr+\int_{0}^{t}
[F_1(r)+F_2(r)X_{{r}^{-}}]\,dN_{r} \ ,\\
X_0=\psi(X_1), \quad t \in [0,1]\ .
\end{array}
\right.
$$     
In each of the following cases, $X$ is a reciprocal process:
\begin{description}
\item [$(1)$]
$\psi$ is constant.
\item [$(2)$]
$\psi (0)=0$, $f_1\equiv F_1\equiv 0$.
\item [$(3)$]
$F_2 \equiv 0$.
\item [$(4)$]
$F_2 > -1$, $f_1=F_1\equiv 0$.
\item [$(5)$]
$F_2 \equiv -1$, $f_1=F_1\equiv 0$. 
\end{description}
\end{thm}
{\em Proof}:
$(1)$ reduces to the case of initial condition, while in $(2)$ the solution is 
identically zero.
Thus in both situations we
obtain a Markov process.
\begin{description}
\item [$(3)$]
In this case the solution is
$$
X_t=A(t)\Big[ X_0+\int_{0}^{t}
\frac {f_1(r)}{A(r)}\,dr+\int_{0}^{t}
\frac {F_1(r)}{A(r)}\,dN_{r}\Big]\ ,
$$
where 
$
A(t)=\exp \{ \int_{0}^{t}f_2(r)\,dr \}
$
as before, and $X_0$ solves
$$
X_0= \psi \Big(A(1)\Big[ X_0+
\int_{0}^{1}
\frac {f_1(r)}{A(r)}\,dr+\int_{0}^{1}
\frac {F_1(r)}{A(r)}\,dN_{r} \Big]\Big)
\ .      
$$
  Defining $Y_t:=\frac{X_t}{A(t)}$, we can write $Y_t=\xi_t+g(\xi_1)$,
  where $g$ is a Borel function and 
$$
  \xi_t=\int_{0}^{1}
  \frac {f_1(r)}{A(r)}\,dr+\int_{0}^{1}
  \frac {F_1(r)}{A(r)}\,dN_{r} 
  \ .
$$
Since $\xi$ has independents increments, 
Lemma \ref{JM} implies that $Y$, and therefore $X$, 
are reciprocal processes. 
\item [$(4)$]
Here the solution is given by
$$
X_t=X_0 \exp \Big\{\int_{0}^{t}f_2(r)\,dr+\int_{0}^{t}
\log (1+F_2(r))\,dN_{r} \Big\}\ ,
$$
where $X_0$ satisfies
$$
X_0=\psi \Big( X_0 \exp \Big\{\int_{0}^{1}f_2(r)\,dr
+\int_{0}^{1}
\log (1+F_2(r))\,dN_{r} \Big\}\Big)\ .
$$
When $\psi (0)=0$, we are in case (1).
If $\psi (0)>0$, then $X_0>0$, and setting
$Y_t:=\log (X_t)$ we obtain
$
Y_t = \log [\psi (e^{\xi_1})] + \xi_t
\ ,
$
where $\xi_t=
\int_{0}^{t} f_2(r)\,dr+\int_{0}^{t}
\log [1+F_2(r)]\,dN_{r}.$
We reach the conclusion as in case (3).
If $\psi (0)<0$, we can proceed analogously.
\item [$(5)$]
In this situation, we obtain
\begin{equation}
\label{case5}
X_t(\omega)=
\left \{
\begin{array}{ll}
\ds
x^{\ast}A(t)\ ,& \mbox{if $S_1(\omega)>1$} \\
\psi (0)A(t) \ind_{[0,S_1(\omega))}(t)\ , &
\mbox{if $S_1(\omega)\leq 1$}\ ,
\end{array}
\right.
\end{equation}
where $x^{\ast}$ solves
$x=\psi (xA(1))$. 
Process (\ref{case5}) can be thought as the solution
  to the initial value problem
$$
  X_t=\eta+\int_0^t f_2(r)X_{r^-}\,dr-\int_0^t X_{r^-}\,dN_r
  \ ,
$$
  where $\eta$ is the random variable 
$$
  \eta(\omega)=
  \left\{
  \begin{array}{ll}
  x^*\ ,& \mbox{if $S_1(\omega)>1$}
  \\
  \psi(0)\ ,& \mbox{if $S_1(\omega)\le 1$} \ .
  \end{array}
  \right.
$$
  Although $\eta$ anticipates the Poisson process, 
  $X$ not only has the reciprocal
  property, but it is in fact a Markov process.
  Indeed, it is immediate to check that $Y_t:=X_t/A(t)$ is
  a Markov chain taking at most three values.
  \qed
\end{description}
\setcounter{equation}{0}
\section{Backward and Skorohod equations}
\label{backsko}
In this Section we consider the backward and Skorohod versions
of our boundary value problems. There are very simple cases 
where the backward equation, even in the initial condition situation,
does not possess a solution. For example, for $k\in\Reals$,  
the equation
$$
  \varphi_t=1+\int_0^t k\varphi_s\,dN_s
  \ ,
$$
leads to 
$\varphi_{S_1}=1+k\varphi_{S_1}$
  at $t=S_1$, 
  which is absurd for $k=1$. 
  In general, for the existence of a solution of 
\begin{equation}
\label{flujobackward}
\varphi_{t}(x)=x+\int_{0}^{t}f(r,\varphi_{r}(x))\,dr
+\int_{0}^{t}F(r,\varphi_{r}(x))\,dN_{r}\ ,
\quad t \in [0,1]\ ,
\end{equation}
it is necessary that the mapping
  $A_r(y):=y-F(r,y)$ be invertible for each $r$. 
\par
Assume now that $f\colon [0,1]\times \Reals \rightarrow \Reals$ satisfies
the hypotheses $(H_1)$ and $(H_2)$ of Section \ref{flows}, 
and that either 
\begin{equation}
\label{Con1}
  \forall t\in[0,1],\ \exists \alpha(t)< 1:\ x>y\Rightarrow F(t,x)-F(t,y)\le \alpha(t)(x-y)
  \ ,
\end{equation}
or
\begin{equation}
\label{Con2}
  \forall t\in[0,1],\ \exists \alpha(t)> 1:\ x>y\Rightarrow F(t,x)-F(t,y)\ge \alpha(t)(x-y)
  \ .
\end{equation}
Then 
there exists a 
unique process 
$\varphi(x)=\{\varphi_t(x), t \in [0,1]\}$
that satisfies the backward equation (\ref{flujobackward}).
This follows from Theorem $5.1$ of Le\'on, Sol\'e and Vives \cite {LSV},
where it is shown that $\varphi$ is a solution to (\ref{flujobackward})
if and only if $\varphi$ is a solution to the forward equation
$$
\varphi_{t}(x)=x+\int_{0}^{t}
f(r,\varphi_{{r}^{-}}(x))\,dr+\int_{0}^{t}
F(r,A_{r}^{-1}(\varphi_{{r}^{-}}(x)))\,dN_{r}\ .
$$
The existence of $A_r^{-1}$
is assured by (\ref{Con1}) or (\ref{Con2}).

We consider now the backward equation with boundary condition
\begin{equation}
\label{D2}
%({\Cal D})
\left \{
\begin{array}{l}
\ds
X_t=X_{0}+\int_{0}^{t}f(r , X_{r})\, d r
+\int_{0}^{t}F(r,X_{r})\,dN_{r} \ ,\\
X_0= \psi (X_1)\ ,\quad t \in [0,1]\ .
\end{array}
\right.
\end{equation} 
\begin{thm}
\label{spD}
Assume that $f$ satisfies hypotheses $(H_1)$ and $(H_2)$ 
of Section \ref{flows}
and that 
\begin{equation}
\label{alphabeta}
  \beta(t)(x-y)\le F(t,x)-F(t,y)\le \alpha(t)(x-y)
  \ ,\quad
  t\in[0,1]
  \ ,\quad
  x>y
  \ ,
\end{equation}
  for some functions $\alpha$ and $\beta$ such that 
$\alpha-1\le \beta\le\alpha<1$.  
Assume 
moreover that $\psi$ satisfies hypothesis $(H_3)$ of Theorem \ref{spB}.
Then (\ref{D2}) admits a unique solution
$X=\{X_t, t \in [0,1]\}$, which is a c\`{a}dl\`{a}g process.
\end{thm}
{\em Proof}:
  By the relation between the forward and the backward
  equation with initial condition given above, the solution to (\ref{D2})
  coincides with the solution to the forward equation with 
  boundary condition
\begin{equation}
%({\Cal D})
\left \{
\begin{array}{l}
\ds
X_t=X_{0}+\int_{0}^{t}f(r , X_{r^-})\, d r
+\int_{0}^{t}F(r,A_r^{-1}(X_{r^-}))\,dN_{r} \ ,\\
X_0= \psi (X_1)\ ,\quad t \in [0,1]\ ,
\end{array}
\right.
\end{equation}
  provided it exists. 
  By Theorem \ref{spB}, it is enough to show 
  there exists a constant $k_2\ge -1$ such that  
  $\tilde F(t,x):=F(t,A_t^{-1}(x))$ satisfies
$\tilde F(t,x)
  -
  \tilde F(t,y)
  \ge 
  k_2
  (x-y)$,
  $\forall
  t\in[0,1]$, 
  $x>y$.
  From (\ref{alphabeta}), 
$$
  0<(1-\alpha(t))(x-y)\le A_t(x)-A_t(y)\le(1-\beta(t))(x-y)
  \ ,
$$
  hence
$$
  \frac{x-y}{1-\beta(t)}
  \le 
  A_t^{-1}(x)-A_t^{-1}(y)
  \le
  \frac{x-y}{1-\alpha(t)}
  \ .
$$
  We find
$$
  \tilde F(t,x)
  -
  \tilde F(t,y)
  \ge 
  \left \{
  \begin{array}{ll}
  \ds
  \frac{\beta(t)}{1-\beta(t)}(x-y)
  \ ,
  &
  \mbox{if $\beta(t)\ge 0$}
  \ ,
  \vspace{2truemm}
  \\
  \ds
  \frac{\beta(t)}{1-\alpha(t)}(x-y)
  \ ,
  &
  \mbox{if $\beta(t)<0$}
  \ ,
  \end{array}
  \right.
$$
  and the conclusion follows.
\qed
\par
  The study of the properties of backward equations can thus be reduced to the
  case of the forward equations when the above condition (\ref{alphabeta}) on 
  $F$ holds true.
  In particular, when $F(t,x)=F_1(t)+F_2(t)x$, we obtain 
$$
  \tilde F(t,x)=\frac{F_1(t)}{1-F_2(t)}
  +\frac{F_2(t)}{1-F_2(t)}\,x
  \ ,
$$
  and condition (\ref{alphabeta}) reduces to $F_2<1$.
\begin{example}
\label{backlinex}
{\sl (Linear backward equation)}.
{\rm
Now consider the problem
$$
%({\Cal LD})
\left \{
\begin{array}{l}
\ds
X_t=X_0+\int_{0}^{t} 
[f_1(r)+f_2(r)X_{{r}}]\,dr +\int_{0}^{t}
[F_1(r)+F_2(r)X_{{r}}]\,dN_{r} \ ,\\
X_0=\psi(X_1)\ , \quad  t \in [0,1]\ ,
\end{array}
\right.
$$
where $f_1,f_2,F_1,F_2\colon[0,1] \rightarrow \Reals$ are
continuous functions with 
$F_2<1$
and $\psi\colon\Reals \rightarrow \Reals$ is a continuous and
non-increasing function. By Theorem \ref{spD},
this problem has unique solution, given by (see Example \ref{OE})
$$
X_t=\eta_{t} \Big[ X_{0}+\int_{0}^{t}
\frac {f_1(r)}{\eta_{r}}\,dr+\int_{0}^{t}
\frac {F_1(r)}
{(1-F_{2}(r))\eta_{r}}\,dN_{r} \Big]\ ,
$$
where                                                
$$
\eta_t=\exp\Big\{ \int_{0}^{t}f_2(r)\,dr
-\int_{0}^{t}\log (1-F_2(r))\,dN_{r} \Big\}\ ,     
$$                  
and $X_{0}$ solves
$$
x= \psi \Big(\eta_{1} \Big[ x+\int_{0}^{1}
\frac {f_1(r)}{\eta_{r}}\,dr+\int_{0}^{1}
\frac {F_1(r)}
{(1-F_{2}(r))\eta_{r}}\,dN_{r} \Big] \Big)\ .
\quad\qed
$$  
}
\end{example}
\par
\bigskip
  We turn now to 
  the Skorohod equation with boundary condition
\begin{equation}
\label{C2}
\left \{
\begin{array}{l}
\ds
X_t=X_{0}+\int_{0}^{t}f(r , X_{r})\, d r
+\int_{0}^{t}F(r,X_{r})\,\delta \tilde N_{r} \ ,\\
X_0= \psi (X_1)\ , \quad t \in [0,1]\ .
\end{array}
\right.
\end{equation} 
We place ourselves in the canonical Poisson space 
$(\Omega,\goth F,P)$ (see e.g. \cite{LSV}, \cite{NV} or \cite{NV2} for a 
more detailed introduction to the analysis in this space).
The 
elements of $\Omega$ are sequences $\omega=(s_1,\dots,s_n)$,
$n\ge 1$,
with $s_j\in [0,1]$, together with a special point $a$. 
The canonical Poisson process is defined in $(\Omega,\goth F,P)$
as the measure-valued process 
$$
  N(\omega)=
  \left\{
  \begin{array}{ll}
  0\ , &\mbox{if $\omega=a$}\ ,
  \\
  \sum_{i=1}^n \delta_{s_i}\ , &\mbox{if $\omega=(s_1,\dots,s_n)$}\ ,
  \end{array}
  \right.
$$
  where $\delta_{s_i}$ means the Dirac measure on $s_i$.
  Any square integrable random variable $H$ in this space can be decomposed 
  in Poisson-It\^o chaos 
  $H=\sum_{n=0}^{\infty} I_n(h_n)$, where $I_n(h_n)$ is the $n$-th multiple
  Poisson-It\^o integral of a symmetric kernel $h_n\in L^2([0,1]^n)$ 
  with respect to the compensated Poisson process $\tilde N$.
  For $u\in L^2(\Omega\times [0,1])$ with decomposition $u_t=\sum_{n=0}^{\infty} I_n(u_n^t)$
  for almost all $t\in[0,1]$,
  Nualart and Vives \cite{NV} 
  define its Skorohod integral as 
  $\delta(u):=\int_0^1 u_s\,\delta\tilde N_s:=
  \sum_{n=0}^{\infty}
  I_{n+1}(\tilde u_n)$,
  where 
  $\tilde u_n$ is the symmetrization of $u_n$ with respect to its $n+1$
  variables,
  provided $u\in\Dom\ \delta$, that means, if
  $\sum_{n=0}^{\infty} (n+1)!\,\|\tilde u_n\|_{L^2([0,1]^{n+1})}^2<\infty$.
\par
  For a process $u$ with integrable paths, define the random variable 
$$
  \phi(u)(\omega):=
  \left\{
  \begin{array}{ll}
  -\int_0^1 u_t(a)\,dt\ , &\mbox{if $\omega=a$}\ ,
  \\
  u_{s_1}(a)-\int_0^1 u_t(s_1)\,dt\ , &\mbox{if $\omega=(s_1)$}\ ,
  \\
  \sum_{j=1}^n u_{s_j}(\hat \omega_j)-\int_0^1 u_t(\omega)\,dt\ , 
    &\mbox{if $\omega=(s_1,\dots,s_n)$, $n>1$}\ ,
  \end{array}
  \right.
$$
  where $\hat \omega_j$ means $(s_1,\dots,s_{j-1},s_{j+1},\dots,s_n)$.
  One can also consider,
  for any random variable $H$ and for almost all $t\in[0,1]$, 
  the random variable
$$
  (\Psi_t H)(\omega)=
  \left\{
  \begin{array}{ll}
  H(t)-H(a)\ , &\mbox{if $\omega=a$}
  \\
  H(s_1,\dots,s_n,t)-H(\omega)\ , &\mbox{if $\omega=(s_1,\dots,s_n)$}\ .
  \end{array}
  \right.
$$
  The following Lemma is shown in 
  Nualart and Vives \cite{NV}. 
\begin{lm}\label{NV-lemma}
  With the notations introduced above, we have 
\begin{description}
\item[(a)]
   If $u\in L^2(\Omega\times[0,1])$, then
  $\phi(u)\in L^2(\Omega)$ if and only if $u\in\Dom\ \delta$, and in that
  case $\delta(u)=\phi(u)$. 
\item[(b)] 
  If $H=\sum_{n=0}^{\infty} I_n(h_n)\in L^2(\Omega)$, 
  then $\Psi H\in L^2(\Omega\times[0,1])$ if and only if
  $\sum_{n=0}^{\infty} n\,n!\,\|h_n\|_{L^2([0,1]^n)}^2<\infty$, and in that case
  $\Psi_t H=\sum_{n=0}^{\infty} (n+1)I_n(h_{n+1}(t,\cdot))$.
\end{description}
\end{lm}
\par
  Two concepts of solution for initial value Skorohod
  equations were introduced in \cite{LCT} by Le\'on, Ruiz de Ch\'avez  
  and Tudor,
  which they called ``strong solution'' and ``$\phi$-solution''. In the 
  latter, the process $F(r,X_r)$ is only required to have square 
  integrable paths, and its integral is interpreted as 
  $\phi(F(r,X_r))$. If $F(r,X_r)$ belongs to $\Dom\ \delta$, 
  Lemma \ref{NV-lemma} (a) ensures that both concepts coincide. 
  We only need here a version of the first notion, which 
  we will call simply ``solution''.
  We supplement the definition in \cite{LCT} with 
  the requirement of c\`adl\`ag
  paths, for the boundary condition to be meaningful.
\begin{de}
{\rm
  A measurable process $X$ is a {\sl solution\/} of
  (\ref{C2}), if 
\begin{description}
\item [$(1)$]
$f(\cdot,X_{\cdot}) \in L^{1}([0,1])$ with probability 1.
\item [$(2)$]
$\ind_{[0,t]}(\cdot)F(\cdot, X_{\cdot}) \in 
\Dom\ \delta $ for almost all $t \in [0,1]$.
\item [$(3)$]
The first equality in (\ref{C2}) is satisfied with 
probability 1, for almost all $ t \in [0,1]$.
\item [$(4)$]
With probability 1, $X$ is c\`adl\`ag and $X_{0}=\psi (X_{1})$.
\qed
\end{description}
}
\end{de}
\begin{thm}
\label{spC}
Let $f,F\colon[0,1]\times \Reals \rightarrow \Reals$
satisfy hypotheses $(H_1)$ and $(H_2)$ of Section \ref{flows} 
with constants $K_1$, 
$M_1$ and $K_2$, $M_2$, respectively. Assume moreover
that $\psi\colon\Reals \rightarrow \Reals$ satisfies
\begin{description}
\item [$(H_3'')$]
$\psi$ is a continuous and bounded function that verifies
one of the following Lipschitz-type conditions:
\begin{description}
\item [$(i)$]
$x>y \Rightarrow \psi (x)-\psi (y)\leq \eta\cdot (x-y)$,
 for some real constant $\eta < e^{-\tilde K}$,
\item [$(ii)$]
$x>y \Rightarrow \psi (x)-\psi (y)\geq \eta\cdot (x-y)$,
for some real constant $\eta > e^{\tilde K}$,
\end{description}
where $\tilde K=K_1+K_2$.
\end{description}
Then (\ref{C2}) admits a unique 
%$\phi$-solution
%$X$, which is also its unique 
solution.
\end{thm}
{\em Proof}:
Under our hypotheses, we can apply Theorems 3.7 and 3.13 of
\cite{LCT} to the equation
\begin{equation}
\label{441}
X_t=\zeta+\int_{0}^{t}f(r,X_{{r}})\,dr+
\int_{0}^{t}F(r,X_{{r}})\,\delta \tilde N_{r} 
\ ,
\end{equation}
where $\zeta$ is a given bounded random variable, and the solution
 has a c\`adl\`ag version
given by
$$
X_{t}=\sum_{n=0}^{\infty}
X_{t}^{n}(\omega)\ind_{[0,1]^{n}}(\omega)\ ,
$$
where
$X^n$ are the respective unique solutions of 
\begin{equation}
\label{442}
X_{t}^{0}(a)=\zeta(a)+\int_{0}^{t}
(f-F)(r,X_{r}^{0}(a))\,dr\ ,
\end{equation}
\begin{equation}
\label{443}
X_{t}^{1}(s_1)=\zeta(s_1)+\int_{0}^{t}
(f-F)(r,X_{r}^{1}(s_1))\,dr+
\ind_{[0,t]}(s_1)F(s_1,X_{s_1}^{0}(a))\ ,
\end{equation}
and for $n \geq 2$, 
\begin{equation}
\label{444}
X_{t}^{n}  (\omega)
=\zeta(\omega)+\int_{0}^{t}
(f-F)(r,X_{r}^{n}(\omega))\,dr
+\sum_{j=1}^{n} \ind_{[0,t]}(s_j)
F(s_j,X_{s_j}^{n-1}(\tilde \omega_{j}))
\ ,
\end{equation}
with
$\omega=(s_1,\dots,s_n)$ and
$\hat \omega_{j}=(s_1,\dots,s_{j-1},s_{j+1},\dots,s_n)$.
\par
Denoting by $X^n(\omega,x)$ the corresponding solutions 
starting at $\zeta\equiv x\in\Reals$,
it is easy to show that 
for any $x_1,x_2 \in \Reals$ with
$x_1 < x_2$,
$$
(x_2- x_1)e^{-\tilde K t} \leq
X_{t}^{n}(\omega,x_2)-
X_{t}^{n}(\omega, x_1) \leq
(x_2 - x_1)e^{\tilde K t}
\ .
$$
These inequalities and $(H_3'')$ imply that
there exists a unique point $x^*$ such that
\begin{equation}
\label{445}
x^*=\psi (X_{1}^{0}(a,x^*))
\ ,
\end{equation}
which we define as $X_0(a)$. Then, $X^0(a,X_0(a))$ satisfies (\ref{442}) 
with $\zeta(a)=X_0(a)$ and 
the boundary condition of (\ref{C2}).
\par
Once we know $X^0$, using equation (\ref{443}) and reasoning 
 similarly, one shows that there exists a unique point $x^*$
such that
\begin{equation}
\label{446}
x^*=\psi (X_{1}^{1}(s_1,x^*))\ ,
\end{equation}
which we define as $X_{0}(s_1)$.
Then, $X^1(s_1,X_0(s_1))$ satisfies (\ref{443}) with $\zeta(s_1)=X_0(s_1)$
and the boundary condition
of (\ref{C2}).
\par
In general, given $\omega=(s_1,\dots,s_n)$, once we know $X^{n-1}$,
and using equation (\ref{444}), one shows that there exists a unique point
$x^*$ such that
$$
x^*=\psi (X_{1}^{n}(\omega,x^*))\ ,
$$
which we define as $X_0(\omega)$.
We have then that $X^n(\omega,X_0(\omega))$ satisfies
(\ref{444}) with $\zeta(\omega)=X_0(\omega)$ and the boundary condition of (\ref{C2}).
\par
Since $\psi$ is bounded, $X_0$ is a bounded random
variable. The process thus constructed clearly satisfies 
(\ref{441}) together with the boundary condition, and the theorem
is proved.
\qed
\begin{re}
\label{B=C}
  {\rm
  Skorohod equations can be converted to forward ones 
  in special situations:
  When $F(t,x)\equiv F(t)$ or when $\psi \equiv x_0\in\Reals$,
  the solution of (\ref{C2}) coincides with the
  solution of
$$
\left \{
\begin{array}{l}
\ds
X_t=X_0+\int_{0}^{t}(f-F)(r,X_{{r}^{-}})\,dr+
\int_{0}^{t}F(r,X_{{r}^{-}})\,d N_{r} \ ,\\
X_0=\psi (X_1)\ ,\quad t \in [0,1]\ .
\quad\qed
\end{array}
\right.
$$ 
}
\end{re}
The results of Sections \ref{eqwbc} and \ref{reciprocal}
are automatically translated to Skorohod equations in the situations of the 
previous remark. In other cases, the inductive construction of the 
solution $X$, in which the value of $X_t(\omega)$ (with $\omega\in[0,1]^n$)
depends on the values $X_t(\omega)$ (with $\omega\in[0,1]^{n-1}$), does not allow 
the equivalence.

In the last five years there have been 
some interest in finding the chaos decomposition of solutions to
several type of equations in Poisson space
(see e.g. \cite{LT}, \cite{LSV}). 
For instance, for 
$$
X_t=x+\int_{0}^{t}f_2(r)X_r\,dr+
\int_{0}^{t}F_2(r)X_r\,\delta \tilde N_{r} 
\ ,
\quad x\in\Reals
\ ,
$$ 
one can find, using Lemma 3.10 of \cite{LT}, the decomposition
$$
X_t=x \exp \Big\{ \int_{0}^{t}(f_2+F_2)(r)\,dr \Big\}
\sum_{n=0}^{\infty}
I_{n}[(\ind_{[0,t]}(\cdot) F_2(\cdot))^{\otimes n} ]/n! \ .
$$
We will give the chaos decomposition of the solution of
two very specific linear equations with boundary conditions. First we discuss 
the resolution of Skorohod linear equations.
\begin{example}
(Linear Skorohod equation).
{\rm
Consider the problem
$$
\label{LC}
\left \{
\begin{array}{l}
\ds
X_t=X_0+\int_{0}^{t} 
[f_1(r)+f_2(r)X_{{r}}]\,dr +\int_{0}^{t}
[F_1(r)+F_2(r)X_{{r}}]\delta \tilde N_{r} \ ,\\
X_0=\psi(X_1)\ , \quad  t \in [0,1]\ ,
\end{array}
\right.
$$
where $f_1,f_2,F_1,F_2\colon[0,1] \rightarrow \Reals$ are
continuous functions, and 
$\psi\colon\Reals \rightarrow \Reals$ satisfies
$(H_3'')$ of the Theorem \ref {spC}. 
\bigskip
To describe $X$, let 
$
Y_{t}(\omega)=\sum_{n=0}^{\infty} Y_{t}^{n}(\omega)
\ind_{[0,1]^{n}}(\omega),
$
the
solution to
$$
Y_t=\int_{0}^{t}[f_1(r)+f_2(r)Y_{r}]\,dr +
\int_{0}^{t}[F_1(r)+F_2(r)Y_{r}]\,\delta \tilde 
N_{r}\ .
$$
Taking into account Remark \ref{B=C},
$Y$ is the solution to the forward equation
$$
Y_t=\int_{0}^{t}[(f_1-F_1)(r)+
(f_2-F_2)(r)Y_{{r}^{-}}]\,dr 
+\int_{0}^{t}[F_1(r)+F_2(r)Y_{{r}^{-}}]\,dN_{r}\ ,
$$
which is given in Example \ref{OE}.
Then $X_t=Y_t+Z_t$, where $Z_t$ satisfies
$$
\left \{
\begin{array}{l}
\ds
Z_t=Z_0+\int_{0}^{t} 
f_2(r)Z_{r}\,dr +\int_{0}^{t}
F_2(r)Z_{r} \delta \tilde N_{r} \ ,\\
Z_0=\psi(Y_1+Z_1)\ , \quad  t \in [0,1]\ .
\end{array}
\right.
$$     
We know (see proof of Theorem \ref{spC}) that $Z_t=\sum_{n=0}^{\infty}Z_{t}^{n}(\omega)
\ind_{[0,1]^{n}}(\omega)$, where, writing 
$\tilde A(t):=\exp\{\int_0^t (f_2-F_2)(r)\,dr\}$,
$$
\frac {Z_{t}^{0}(a)}{\tilde A(t)}= Z_0(a)
$$
and $Z_0(a)$ is the solution to
$$
x =\psi (\tilde A(1)x + Y_1(a))\ .
$$
For $\omega =(s_1)$,
$$
\frac {Z_{t}^{1}(s_1)}{\tilde A(t)}=
Z_{0}(s_1)\ind_{[0,s_1)}(t)+
[Z_{0}(s_1)+F_{2}(s_1)Z_{0}(a) ]\ind_{[s_1,1]}(t)
$$
and $Z_{0}(s_1)$ solves
$$
x=\psi \big( \tilde A(1)[x+F_{2}(s_1))Z_{0}(a)]+Y_{1}(s_1)\big)\ .
$$
In general, for $\omega=(s_1,\dots ,s_n)$ with
$0<s_1<\dots <s_n<1$ we have
\begin{align*}
\frac {Z_{t}^{n}(\omega)}{\tilde A(t)}=
{}&
Z_{0}(\omega)\ind_{[0,s_1)}(t)+\dots 
\\
&\dots+\Big[ Z_{0}(\omega) +
\sum_{k=1}^{i}
\sum_{\scriptstyle
j_{1},\dots,j_{k}=1 \atop \scriptstyle
\text{distinct}
}^{i}F_{2}(s_{j_1})\dots F_{2}(s_{j_k})
Z_{0}(\tilde  \omega_{(s_{j_1},\dots ,s_{j_k})})\Big ]
\ind_{[s_i,s_{i+1})}(t)+\dots 
\\
&
\dots +\Big[ Z_{0}(\omega) +
\sum_{k=1}^{n}
\sum_{\scriptstyle
j_1,\dots ,j_k=1 \atop \scriptstyle
\text{distinct}
} ^{n}F_{2}(s_{j_1})\dots F_{2}(s_{j_k})
Z_{0}(\tilde  \omega_{(s_{j_1},\dots ,s_{j_k})})\Big ]
\ind_{[s_n,1]}(t)\ ,
\end{align*} 
where
$$
\hat \omega_{(s_{j_1},\dots ,s_{j_k})}=
(\dots ,s_{j_{1}-1},s_{j_{1}+1},\dots ,s_{j_{k}-1},s_{j_{k}+1},\dots ) \ ,
\quad \hat \omega_{(s_1,\dots ,s_n)}=a\ ,
$$
and $Z_{0}(\omega)$ is the solution to
$$
x=\psi \Big( \tilde A(1)\Big[x +
\sum_{k=1}^{n}
\sum_{\scriptstyle
j_1,\dots ,j_k=1 \atop \scriptstyle
\text{distinct}
} ^{n}F_{2}(s_{j_1})\dots F_{2}(s_{j_k})
Z_{0}(\tilde  \omega_{(s_{j_1},\dots ,s_{j_k})})\Big]
+ Y_{1}(\omega) \Big ) \ .
\quad\qed
$$ 
}
\end{example}
\begin{example}
(Chaos decompositions)
{\rm
\nopagebreak
\begin{description}
\item [$(1)$]
Consider the problem
$$
\left \{
\begin{array}{l}
\ds
X_t=X_0+\int_{0}^{t}[f_1(r)+f_2(r)X_{{r}}]\,dr
+\int_{0}^{t}F_1(r)\,\delta \tilde N_{r} \ ,\\
X_0=aX_1+b\ , \quad t \in [0,1]\ ,
\end{array}
\right.
$$  
where $f_1,f_2,F_1$
are continuous functions, and $a,b \in \Reals$ with
$a \neq \exp \{ \int_{0}^{1} f_2(r)\,dr \}$.
Its solution coincides with the solution of the forward 
equation
$$
\left \{
\begin{array}{l}
\ds
X_t=X_0+\int_{0}^{t}[(f_1-F_1)(r)+f_2(r)X_{{r}^{-}}]\,dr
+\int_{0}^{t}F_1(r)\,dN_{r} \ ,\\
X_0=aX_1+b\ , \quad t \in [0,1]\ ,
\end{array}
\right.
$$
which is 
$$
X_t=A(t)\Big[X_0 + \int_{0}^{t}\frac {(f_1-F_1)(r)}{A(r)}\,dr
+ \int_{0}^{t}
\frac {F_1(r)}{A(r)}\,dN_{r} \Big]\ ,
$$
where $A(t)=\exp \{\int_{0}^{t}f_2(r)\,dr \}$ and
$$
X_0=\frac {aA(1)}{1-aA(1)} \Big[\int_{0}^{1}
\frac {(f_1-F_1)(r)}{A(r)}\,dr+ 
\int_{0}^{1}
\frac {F_1(r)}{A(r)}\,dN_{r}\Big]
+\frac {b}{1-aA(1)}\ .
$$       
Therefore $X_t$ belongs the first order chaos and 
\begin{align*}
X_t={}& \frac {bA(t)}{ 1-aA(1)}+
\int_{0}^{1}A(t) \Big(\ind_{[0,t]}(r)+
\frac {aA(1)}{1-aA(1)}\Big)
\frac {f_1(r)}{A(r)} \,dr
\\
&+I_{1} \Big[ A(t)\Big( \ind_{[0,t]}(\cdot)+
\frac {aA(1)}{1-aA(1)} \Big)\frac {F_1(\cdot)}{A(\cdot)}
\Big ]\ .
\end{align*}
\item [$(2)$]
Consider the problem
$$
\left \{
\begin{array}{l}
\ds
X_t=X_0+\int_{0}^{t}f_2(r)X_{{r}^{-}}\,dr
-\int_{0}^{t}X_{{r}^{-}}\,dN_{r} \ ,\\
X_0= \psi(X_1)\ , \quad t \in [0,1]\ ,
\end{array}
\right.
$$
where $f_2\colon[0,1] \rightarrow \Reals$
is a continuous function, and
$\psi\colon\Reals \rightarrow \Reals $ is a continuous and
non-increasing function. The solution is
$$
X_t(\omega)=
\left \{
\begin{array}{ll}
\ds
x^{\ast}A(t)\ , & \mbox{if $S_1(\omega)>1$} \\
\psi (0)A(t) \ind_{[0,S_1(\omega))}(t)\ ,
&\mbox{if $S_1(\omega)\leq 1$}\ ,
\end{array}
\right.
$$
where $x^{\ast}$ is the unique solution to
$x=\psi (xA(1))$. We can write in Poisson space
$$
X_t=A(t)\Big[ x^{\ast} \ind_{\{a\}}+
\psi (0) \ind_{\{t<S_1 \}}
\Big]\ .
$$
Using Lemma \ref{NV-lemma} (b) one obtains the following chaos decomposition:
\begin{align*}
X_0={} &\psi (0)+(x^{\ast}-\psi(0))e^{-1}
\sum_{n=0}^{\infty} \frac{(-1)^{n}}{n!}I_{n}(1) \ ,
\\
X_1={} &A(1) x^{\ast}e^{-1}
\sum_{n=0}^{\infty} \frac{(-1)^{n}}{n!}I_{n}(1) \ ,
\end{align*}
and for $t \in (0,1)$ we have
$$
X_t= A(t)\Big\{x^{\ast}e^{-1}
\sum_{n=0}^{\infty} \frac{(-1)^{n}}{n!}I_{n}(1) +
\psi(0)e^{-1}\sum_{n=0}^{\infty}\frac{(-1)^{n} 
}{n!}I_{n} \big[(\ind_{[0,t]}(\cdot))^{\otimes n} \big] \Big\}\ .
\quad\qed
$$ 
\end{description}
}
\end{example}

\end{document}